# RESEARCH ANNOUNCEMENT

# BIZARRE TOPOLOGY IS NATURAL IN DYNAMICAL SYSTEMS

JUDY A. KENNEDY AND JAMES A. YORKE

ABSTRACT. We describe an example of a $C^\infty$ diffeomorphism on a 7–manifold which has a compact invariant set such that uncountably many of its connected components are pseudocircles. (Any 7–manifold will suffice.) Furthermore, any diffeomorphism which is sufficiently close (in the $C^1$ metric) to the constructed map has a similar invariant set, and the dynamics of the map on the invariant set are chaotic.

Compact invariant sets of typical diffeomorphisms are often pictured as being points or manifolds or as sets which are locally the product of a manifold and a Cantor set. Such views would mislead one to believe more exotic topological structures are either exceptional or absent. This is not the case. We report here on an example of a $C^\infty$ diffeomorphism $F$ with the following features:

(1) There is a compact invariant set $\mathbb{B}$ (that is, $F(\mathbb{B}) = \mathbb{B}$) which is maximal in the sense that there is a neighborhood $U$ of $\mathbb{B}$ such that if $F^n(x)$ is in $U$ for all integers $n$, then $x$ is in $\mathbb{B}$.
(2) There is a point in $\mathbb{B}$ whose trajectory is dense in $\mathbb{B}$ under $F$.
(3) The number of connected components of $\mathbb{B}$ is uncountable. Most of the components in the sense of category are, from most perspectives, topologically bizarre objects known as pseudocircles.
(4) For each diffeomorphism $G$ in some $C^1$ neighborhood of $F$, the maximal invariant set $\mathbb{B}_G$ of some neighborhood $U_G$ has a similar topological structure, in that it consists of an uncountable collection of connected components, and most of the components of the invariant set are pseudocircles.

What then is a pseudocircle, and what is bizarre about its topology? Pseudocircles have topological dimension 1, but they are arcwise totally disconnected: every continuous map of an interval into a pseudocircle is constant. A **continuum** is a compact, connected subset of a metric space. A **subcontinuum** of a continuum

Received by the editors August 24, 1993, and, in revised form, November 23, 1994.
1991 *Mathematics Subject Classification*. Primary 54F20, 54F50, 54F12.
This research was supported in part by NSF grant DMS–9006931 and AFOSR grant F49620–J–0033.

is a subset of the continuum which is itself a continuum, i.e., it is connected and compact. A continuum is **decomposable** if it is the union of two of its proper subcontinua. Otherwise the continuum is **indecomposable**. For example, the interval [0,1] is a decomposable continuum because it is $[0, 2/3] \cup [1/3, 1]$, and hence it is the union of two proper subcontinua of [0,1]. In fact, obvious examples of decomposable continua include the connected, compact manifolds. Examples of indecomposable continua are perhaps not so obvious, but these continua not only exist but are **most** continua (in the sense of category) [B] and arise quite frequently in "nice" dynamical systems as invariant sets.

The first example of an indecomposable continuum is due to L. E. J. Brouwer, who concocted it in 1910 as a strange example to disprove a conjecture of Schoenflies that the common boundary of two disjoint regions in the plane had to be decomposable. By the 1920s these objects were being studied, not as examples of just "interesting pathology", but because they have interesting topological properties. In 1932, an indecomposable continuum arose as an invariant set in a dynamical system. G. D. Birkhoff's "remarkable curve" is the invariant set for a diffeomorphism on an annulus, and Marie Charpentier later proved that this curve is an indecomposable continuum. (See [Bi] and [C].) Since the 1960s these continua have arisen in many contexts in many dynamical systems. Examples include the global attractor for the Smale horseshoe map (one of the main motivating examples of modern dynamics), the solenoids (these arise in connection with differential equations and flows), many of the inverse limit systems used to model some dynamical systems (a tool introduced by R. Williams ([W1] and [W2])), and invariant sets associated with the forced van der Pol equations [BG] and forced damped pendulum equations [KY1]. Indecomposable continua arise under mild conditions, when the stable and unstable manifolds of a hyperbolic fixed point $p$ intersect at a point $q \neq p$ (a theorem due to Marcy Barge [Ba], from which it follows that at many parameter values the Hénon and Ikeda attractors are indecomposable). Hence, what Brouwer introduced as pathological is in fact now recognized as mundane in dynamics. Our goal is to carry this process one step further.

All the examples of indecomposable continua mentioned thus far have in common that they all contain arcs, i.e., continua homeomorphic to unit intervals, and they all therefore contain decomposable continua as proper, nowhere dense subcontinua. By 1920, Knaster had described an example of a hereditarily indecomposable continuum. A **hereditarily indecomposable** continuum is a continuum with the property that each of its subcontinua is itself an indecomposable continuum. Therefore, none of the indecomposable continua discussed so far in connection with dynamics is hereditarily indecomposable. However, pseudocircles are hereditarily indecomposable. In 1982, M. Handel [H] gave an example of a $C^\infty$ diffeomorphism of the plane with a strange attractor that is a pseudocircle. Several modifications of Handel's example have appeared since, perhaps the most notable of which is due to M. Herman [He] and is a $C^\infty$ diffeomorphism $f$ with an invariant pseudocircle $\Gamma$ that divides the plane into two regions such that $f$ is analytic on the bounded component of $\mathbb{R}^2 - \Gamma$. Here again, the object was to demonstrate how complicated the topology of an invariant set can be, even in the presence of strong derivative conditions. These constructions appear to be exceptional in that if these examples are perturbed slightly, the result may well be a topologically much nicer invariant set. Perturbations, even $C^\infty$ ones, change the topology of the attractor dramatically. Our goals have been different, in that we can show that this exotic topology not only occurs, it occurs even for an open set of "nice" dynamical systems. Our example is a $C^\infty$ diffeomorphism on a 7–manifold (any 7–manifold will suffice) and has an invariant set that contains many connected components that are pseudocircles and is perturbable in the sense that for $C^1$ maps sufficiently close to the example map, the invariant set again has uncountably many connected components, most of which are pseudocircles. (See [KY3].) In our example it is important to note that unlike the examples of Handel and Herman, the pseudocircles are not invariant, although the union of all the pseudocircles in the invariant set is itself invariant.

What does a hereditarily indecomposable continuum look like? How can they be recognized? An indecomposable continuum can be alternatively characterized as a continuum with the property that each of its proper subcontinua is nowhere dense in the continuum. From this it follows that if $K$ is an indecomposable continuum in a manifold and $N$ is a closed neighborhood of a point of $K$ and $N$ intersects $K$ but does not contain $K$, then $K \cap N$ necessarily has uncountably many components. If $K$ is hereditarily indecomposable, then each of the components of $K \cap N$ is indecomposable, and if $N'$ is a closed neighborhood contained in $N$ and $L$ is a component of $K \cap N$ that intersects the interior of $N'$ but is not contained in $N'$, then $L \cap N'$ has uncountably many connected components. This procedure can be repeated ad infinitum. A hereditarily indecomposable continuum $K$ has the extraordinary property that if $L$ and $M$ are subcontinua of $K$ and $L \cap M \neq \varnothing$, then either $L \subseteq M$ or $M \subseteq L$. R. H. Bing [B] proved that most continua (in the sense of category) are hereditarily indecomposable. Indecomposable continua generally arise in dynamical systems as a result of stretch–contract–fold type behavior in the map. In order to get hereditarily indecomposable continua, layers and layers of folds, of all sizes, must be formed—to the point that it is not possible to change the topology by doing more folding.

A continuum $X$ is **circlelike** if for each $\epsilon > 0$, there exists a continuous map $f_\epsilon$ from $X$ to the unit circle $S^1$ such that for each point $z$ in $S^1$, the set $f_\epsilon^{-1}(z)$ has diameter less than $\epsilon$. Pseudocircles can be defined as hereditarily indecomposable, circlelike continua. Note that this means that it is possible to approximate these continua, in a strong way, by simple closed curves. (See [KY2] or [KY3] for more details.) Since pseudocircles have topological dimension 1, each is homeomorphic to a subset of $\mathbb{R}^3$. All of the pseudocircles that occur as components of the invariant set of our diffeomorphism on a 7–manifold are topologically equivalent to each other, but none is homeomorphic to a subset of the plane. (Some pseudocircles are homeomorphic to subsets of the plane, and all such pseudocircles are homeomorphic. There are uncountably many topologically distinct pseudocircles in $\mathbb{R}^3$.)

Pseudocircles and another hereditarily indecomposable continuum called a pseudoarc have been studied widely. To someone with the usual Euclidean intuition, hereditarily indecomposable continua seem to have no structure. Quite the reverse is true: although a totally different intuition is required, these objects have rich structure and many interesting properties. For more information and references regarding pseudocircles and other hereditarily indecomposable continua we refer the reader to [B], [KY2], and [KY3].

The starting point of all our constructions was our discovery that a pseudocircle could be constructed using only two relatively simple maps. A pseudocircle is generally constructed by taking a countable intersection of nested annuli $A_k$, $k = 1, 2, \ldots$, such that the embedding of $A_{k+1}$ in $A_k$ is severely "crooked" (a term that is defined later). Suppose that we define the annulus $A_1$ to be $S^1 \times [-2, 2]$. In the simplest version of our construction, we find an infinite sequence of maps $\{f_i\}$ on $S^1 \times \mathbb{R}^1$ such that

$$\bigcap_{n=1}^{\infty} (f_n \cdot f_{n-1} \cdot \cdots \cdot f_1)^{-1}(A_1)$$

is a nested sequence of annuli whose intersection is a pseudocircle. We write

$$A_n = (f_n \cdot f_{n-1} \cdot \cdots \cdot f_1)^{-1}(A_1).$$

Next we describe a property of the sequence $\{A_k\}$ which guarantees that $\bigcap A_k$ is a pseudocircle. It is perhaps amazing that the following complicated construction can be achieved using only two relatively simple smooth maps.

Let $A_k$, for $k = 1, 2, \ldots$, be a nested sequence of annuli such that each $A_k$ is partitioned into $N_k$ "rectangles" $R_{k,i}$ (for $i$ an element of $Z/N_k$), that is, closed, connected regions such that $R_{k,i}$ intersects $R_{k,j}$ for $i < j$ if and only if $i + 1 = j$, or $i = 0$ and $j = N_k - 1$. We further assume that

$$\lim_{k \to \infty} \max_{j} \text{ diameter } (R_{k,j}) = 0.$$

Let $\tilde{A}_1$ be a noncompact covering space of $A_1$ and, for $k > 1$, let $\tilde{A}_k$ be the lift of $A_k$ and let $\tilde{R}_{k,j}$ be a lift of $R_{k,j}$ (for $j$ in $Z/N_k$). The crucial questions here are (1) How does $A_{k+1}$ lie in $A_k$? and (2) How can such a construction appear in dynamical processes?

We say that $S$ is a **segment** of $\tilde{A}_k(A_k)$ if it is a compact, connected union of rectangles in $\tilde{A}_k(A_k)$. We say the segment $S$ of $\tilde{A}_{k+1}$ has a $(j_0, j_1)$ **wiggle** (where $j_1 > j_0 + 2$) if $S$ intersects $\tilde{R}_{k,j_0}$ and $\tilde{R}_{k,j_1}$ and if every path that lies wholly in $S$ and starts in $\tilde{R}_{k,j_0}$ and ends in $\tilde{R}_{k,j_1}$ has the following property. In traveling from $\tilde{R}_{k,j_0}$, it must pass through $\tilde{R}_{k,j_1-1}$ and then back to $\tilde{R}_{k,j_0+1}$ before entering $R_{k,j_1}$. We say a segment $S$ of $\tilde{A}_{k+1}$ is **crooked** in $\tilde{A}_k$ if whenever $S$ intersects $\tilde{R}_{k,j_0}$ and $\tilde{R}_{k,j_1}$, for $j_1 > j_0 + 2$, then $S$ has a $(j_0, j_1)$ wiggle. We say that $A_{k+1}$ is **crooked** in $A_k$ if the lift $\tilde{S}$ of each segment $S$ of $A_{k+1}$ is crooked in $\tilde{A}_k$. In other words, if $A_{k+1}$ is crooked in $A_k$, then $A_{k+1}$ must reverse directions in $A_k$, and, the more rectangles $A_k$ has, the more times $A_{k+1}$ must reverse directions in $A_k$. In fact, adding just one rectangle to $A_k$ nearly triples the number of direction reversals.

If the annuli $A_k$ satisfy all these properties we have described, then we are guaranteed that $\bigcap A_k$ is a pseudocircle. (See [KY2], [KY3], and [B].) Having seen how it is possible to have a sequence of nested annuli $\{A_k\}$ constructed so that $\bigcap A_k$ is a pseudocircle, we need to know how such a construction appears in our dynamical processes.

Dynamical systems often exhibit expansion or contraction in one or more directions, as well as shearing and rippling behaviors, and attracting or repelling fixed points. Our diffeomorphism is constructed from a finite and explicit set of these behaviors. The maps for our simplest construction are defined on the cylinder $S^1 \times \mathbb{R}^1$, and are called $T$ (for two, because this map wraps twice) and $W$ (for wiggle, because a wiggle is what $W$ gives us). Specifically, we are interested in the behavior of the maps on the annulus $A = S^1 \times [-2, 2]$. Define the "two" map by

$$T(x, y) = (2x \mod 2\pi, 8y).$$

To define $W$, let $M \geq 2^9$ and define the map $s$, which has an invariant sine curve and is expanding in the $y$ direction, by

$$s(x,\, y) = (x,\, My - (M-1)\sin(x)),$$

and the shear map $\sigma$ by

$$\sigma(x,\, y) = (x - 2\pi y,\, y),$$

and, finally,

$$W(x,\, y) = s \circ \sigma(x,\, y).$$

Note that $T$ and $W$ are analytic maps, but while $W$ is a diffeomorphism, $T$ is exactly two–to–one and only locally a diffeomorphism. We can apply suitably modified versions (no modular arithmetic) of both $T$ and $W$ to the lift $\tilde{A}_1$ of $A$, as well as to $A$. (Since this results in no ambiguity, we call both the maps and their lifts by the same names, $T$ and $W$.)

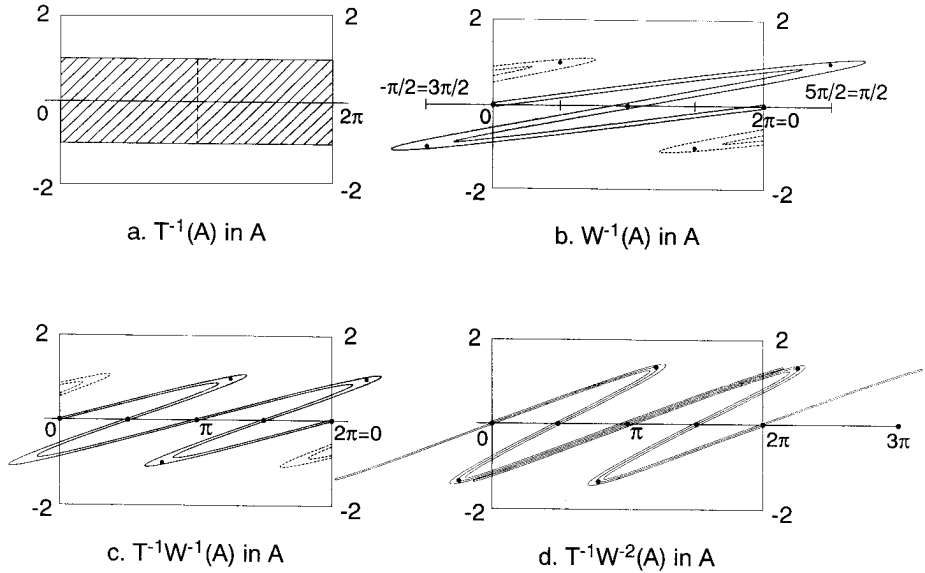

Figure 1

Since $A \subset$ interior $(T(A))$ and $A \subset$ interior $(W(A))$, $T^{-1}(A)$ is a continuum contained in the interior of $A$, as is $W^{-1}(A)$. (See Figure 1.) Suppose that for each $i$, $m_i$ and $n_i$ are positive integers. Define $f_i = W^{n_i} \circ T^{m_i}$. Thus, for each $i$, $f_i^{-1}(A) = T^{-m_i} \circ W^{-n_i}(A)$ is a continuum in the interior of $A$, too. Applying $W^{-1}$ to $A$ (or $\tilde{A}_1$) turns the annulus into a large wiggle—its role is to introduce large folds. If we were to continue applying $W^{-1}$ to $A$ over and over again, the result would be more and more folding, i.e., wiggles inside of wiggles inside of wiggles. But the wiggles stay large in a sense: the only direction reversals occur near the fixed points $p_0 = (\pi/2,\, 1)$ and $p_1 = (-\pi/2,\, -1)$. The result is that $\bigcap_{n \geq 0} W^{-n}(A)$ is an indecomposable continuum. However, it is not hereditarily indecomposable, because we do not have crookedness. If we were to partition $\tilde{A}_1$ into (generalized) rectangles and then partition $W^{-n}(\tilde{A}_1)$ into (generalized) rectangles, we would have only $(j_0,\, j_1)$ wiggles for $j_0$ corresponding to the rectangle containing $p_0$, and $j_1$ corresponding to the rectangle containing $p_1$.

This, of course, is where $T$, or more precisely, $T^{-1}$ comes in. Applying $T^{-1}$ to, say, $W^{-1}(A)$ takes our large wiggle and turns it into exactly two smaller wiggles. (See Figure 1.) Applying $T^{-1}$ to $W^{-2}(A)$ yields two smaller images of the original, too, but this time we get "wiggling" between $-\pi/2$ and $\pi/2$, as well as between $-\pi/4$ and $\pi/4$ and between $3\pi/4$ and $5\pi/4$. Therefore, we can use $T^{-1}$ to adjust wiggle sizes and, more importantly, to get crookedness. (The details of this intuitive discussion can be found in [KY2]. This discussion is not intended to give a rigorous proof, but rather to help the reader understand the ideas of the construction.)

By applying increasingly larger blocks of $T$'s and $W$'s, we can therefore guarantee that for $f_i = W^{n_i} \circ T^{m_i}$, $A_{i+1} = f_i^{-1}(A_i)$ is crooked in $A_i$. The end result is a pseudocircle. Having obtained this pseudocircle, it is not difficult to believe that many different choices of sequences of $T$'s and $W$'s also yield pseudocircles. In fact, in some sense, "most" choices of sequences of $T$'s and $W$'s yield pseudocircles. But, how do we make this precise, and how do we turn all this into one dynamical system, rather than just sequences of maps?

Our next step is to add another factor space, another copy of $S^1$, and a map $g$ on that space. Define $g(z) = 3z \mod 2\pi$ for $z$ in $S^1$. Then choose two disjoint intervals $I_1 = [0, \pi/2]$ and $I_2 = [\pi, 3\pi/2]$ in $S^1$. Note that $g(I_1)$ contains $I_1 \cup I_2$ and is $g|I_1$, one–to–one, and $g(I_2)$ and $g|I_2$ have similar properties. Let

$$C = \{z \in S^1 | g^n(z) \text{ is in } I_1 \cup I_2 \text{ for each } n \geq 0\}.$$

It follows from standard arguments that $C$ is a Cantor set in $S^1$, and, for each sequence $\alpha = \{\alpha_0, \alpha_1, \alpha_2, \dots\}$ consisting of 1's and 2's, there corresponds a unique point $z$ in $C$ such that $\{z\} = \bigcap_0^\infty g^n(I_{\alpha_n})$.

We are finally ready to define our dynamical system: For $(x, y, z)$ in $S^1 \times \mathbb{R}^1 \times (I_1 \cup I_2)$, define

$$F(x, y, z) = \begin{cases} (T(x, y), g(z)) & \text{if } z \text{ is in } I_1; \\ (W(x, y), g(z)) & \text{if } z \text{ is in } I_2. \end{cases}$$

The map is not specified outside $S^1 \times \mathbb{R}^1 \times (I_1 \cup I_2)$, but can be appropriately extended to any 3–manifold in which this manifold is embedded. (See Figure 2.) The third factor $S^1$ plays the role of a function selector, and either $T$ or $W$ is applied to the first and second coordinates of a point depending on whether the third coordinate of the point is in $I_1$ or $I_2$. Thus, the invariant set is

$$\mathbb{B} = \bigcap_0^\infty F^n(S^1 \times \mathbb{R}^1 \times (I_1 \times I_2));$$

the projection $\pi_3(\mathbb{B})$ of $\mathbb{B}$ to the third factor is the Cantor set $C$; and if $z$ is in $C$ and $\alpha$ is the sequence of 1's and 2's associated with $z$, then $B_z = \bigcap_0^\infty F^n(S^1 \times \mathbb{R}^1 \times I_{\alpha_n})$ is a nonempty, connected component of $\mathbb{B}$. We refer to the $B_z$'s as the slices of $\mathbb{B}$. Note that $\pi_3(B_z) = \{z\}$. Thus, we can think of each slice $B_z$ of $\mathbb{B}$ as being a subset of $A$, since choosing a point $z$ in $C$ is equivalent to choosing a sequence of 1's and 2's and therefore a sequence of maps $T$ and $W$ to apply to $A$. Some of the slices $B_z$ are pseudocircles, and some are simpler continua. However, we prove in [KY2] that for a dense $G_\delta$–set (or a residual set) of slices of $\mathbb{B}$, those slices are pseudocircles.

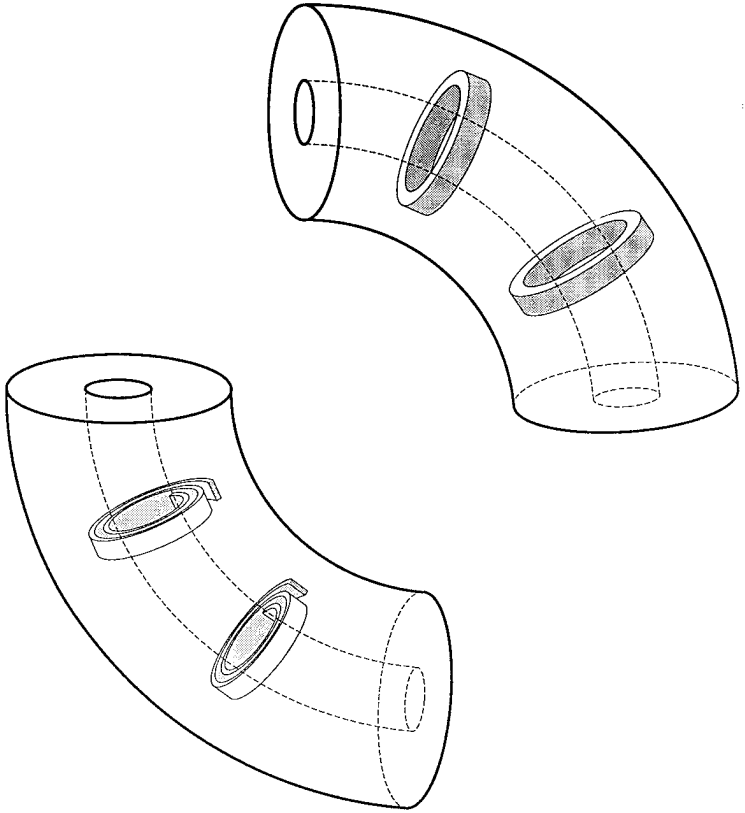

Figure 2

The construction we have just described is a simplified version of the construction given in [KY2]. The $C^\infty$ dynamical system here, unlike that in [KY2], is not perturbable, in the sense that there are maps arbitrarily $C^1$ close to $F$ which do not have a similar repelling invariant set. It is possible for perturbed versions of $F$ to have much simpler, topologically, repelling invariant sets. They may well not have any pseudocircle components. Although the $C^\infty$ map constructed in [KY2] does have the property that $C^1$ maps sufficiently close to that map have a repelling invariant set $\mathbb{B}$ with uncountably many connected components, and most of those components are pseudocircles, that map is not a diffeomorphism because neither $T$ nor $g$ is one–to–one.

In order to construct a $C^\infty$ diffeomorphism, we replace the two circle factors with solid tori. This allows enough "room" to retain the information lost in the 3–manifold construction and gives us a construction on a 7–manifold. Since in the 7–manifold construction of the diffeomorphism (call it $F$), both expansion and contraction occur (in different directions), the resulting invariant set $\mathbb{B}$ is the set of all points that remain in the manifold $D \times [-2, 2] \times D$, where $D$ denotes the solid torus, for all time (both positive and negative). We do not go further into the details of that construction here, as they are considerable, but refer the reader to [KY3]. We do note that, as before, the diffeomorphism $F$ has an invariant set $\mathbb{B}$ possessing the three properties listed at the beginning of the paper. The proof that $\mathbb{B}$ contains a point with a dense trajectory (in $\mathbb{B}$) can be found in [KY4], as well as further discussion of the dynamics of $F$ on $\mathbb{B}$.

## References

[BG]  M. Barge and R. Gillette, *Indecomposability and dynamics on invariant plane separating continua*, Contemp. Math. **117** (1191), 13–38.

[Ba]  M. Barge, *Homoclinic intersections and indecomposability*, Proc. Amer. Math. Soc. **101** (1987), 541–544.

[B]  R. H. Bing, *Concerning hereditarily indecomposable continua*, Pacific J. Math. **1** (1951), 43–51.

[Bi]  G. D. Birkhoff, *Sur quelques courbes fermées remarquables*, Bull. Math. Soc. France **60** (1932), 1–28.

[C]  M. Charpentier, *Sur quelques propriétés des de M. Birkhoff*, Bull. Math. Soc. France **62** (1934), 193–224.

[H]  M. Handel, *A pathological area preserving $C^\infty$ diffeomorphism of the plane*, Proc. Amer. Math. Soc. **86** (1982), 163–168.

[He]  M. Herman, *Construction of some curious diffeomorphisms of the Riemann sphere*, J. London Math. Soc. **34** (1986), 375–384.

[KY1]  J. Kennedy and J. A. Yorke, *Basins of Wada*, Phys. D **51** (1991), 213–225.

[KY2]  ______, *Pseudocircles in dynamical systems*, Trans. Amer. Math. Soc. **343** (1994), 349–366.

[KY3]  ______, *Pseudocircles, diffeomorphisms, and perturbable dynamical systems*, Ergodic Theory Dynamical Systems (to appear).

[KY4]  ______, *The dynamics of some diffeomorphisms which admit complicated invariant sets*, in preparation.

[W1]  R. Williams, *One-dimensional wandering sets*, Topology **6** (1967), 473–487.

[W2]  ______, *Expanding attractors*, Inst. Hautes Études Sci. Publ. Math IHES **43** (1974), 169–203.

Department of Mathematical Sciences, University of Delaware, Newark, Delaware 19716
   *E-mail address*: `jkennedy@brahms.udel.edu`

Institute for Physical Science and Technology, University of Maryland, College Park, Maryland 20742
   *E-mail address*: `yorke@glue.umd.edu`